\documentclass[letterpaper, 10 pt, conference]{ieeeconf}  

\IEEEoverridecommandlockouts                              

\usepackage[utf8x]{inputenc}

\usepackage[greek,english]{babel}
\usepackage{gfsporson}

\usepackage{balance}
\usepackage{graphicx}
\usepackage{siunitx}

\usepackage{pgfplots}
\usepackage{tikz}
  \pgfplotsset{compat=newest}
  \usetikzlibrary{plotmarks}
  \usetikzlibrary{arrows.meta}
  \usepgfplotslibrary{patchplots}
  \usepgfplotslibrary{groupplots}

\newlength\figH
\newlength\figW
\usepackage{amsmath} %
\usepackage{amsfonts}
\usepackage{amssymb}  %
\usepackage{mathrsfs}
\usepackage{url}

\makeatletter
\let\NAT@parse\undefined
\makeatother

\let\ieeebibliography\thebibliography

\usepackage[numbers,sort&compress]{natbib}
\usepackage[colorlinks]{hyperref}
\usepackage[commentmarkup=uwave,final]{changes}
\definechangesauthor[color=purple]{TC}
\definechangesauthor[color=violet]{JO}
\definechangesauthor[color=blue]{JO2}
\renewcommand\thebibliography[1]{\ieeebibliography{#1}}

\usepackage[IEEE]{altthm-styles}

\DeclareMathOperator*{\argmin}{arg\,min}

\title{\LARGE \bf
Ca\textporson{Σ}oS: A nonlinear sum-of-squares optimization {suite}$^\star$
}

\usetikzlibrary{calc}
\newcommand*\circled[1]{\tikz[baseline=(char.base)]{
    \node[shape=circle, draw, inner sep=1pt, 
        minimum height={\f@size*1.6},] (char) {\vphantom{WAH1g}#1};}}

\author{Torbjørn Cunis and Jan Olucak
\thanks{$^\star$This work was not supported by any organization \comment[id=JO]{Maybe add DASO/QUASAR/ATLAS? At least partially?}}
\thanks{The authors are with the Institute of Flight Mechanics and Control, University of Stuttgart, 70569 Stuttgart, Germany. 
        \url{tcunis@ifr.uni-stuttgart.de}, \url{jan.olucak@ifr.uni-stuttgart.de}}%
}

\begin{document}

\maketitle
\thispagestyle{empty}
\pagestyle{empty}

\begin{abstract}
\let\textgreek\textporson
We present Ca\textgreek{Σ}oS, the first MATLAB software specifically designed for nonlinear sum-of-squares optimization. Ca\textgreek{Σ}oS' symbolic polynomial algebra system allows to formulate parametrized sum-of-squares optimization problems and facilitates their fast, repeated evaluations. To that extent, we make use of CasADi's symbolic framework and realize concepts of monomial sparsity, linear operators (including duals), and functions between polynomials. Ca\textgreek{Σ}oS provides interfaces to the conic solvers SeDuMi, Mosek, and SCS as well as methods to solve quasiconvex optimization problems (via bisection) and nonconvex optimization problems (via sequential convexification). Numerical examples for benchmark problems including region-of-attraction and reachable set estimation for nonlinear dynamic systems demonstrate significant improvements in computation time {compared to existing toolboxes.} Ca\textgreek{Σ}oS is available open-source at \url{https://github.com/ifr-acso/casos}.
\end{abstract}

\section{Introduction}
Polynomial optimization subject to the cone of sum-of-squares (SOS) polynomials \cite{parillo2003} provides a powerful tool to tackle various problems that arise in nonlinear systems and control theory. Such problems include region-of-attraction estimation \cite{chakrabortyEtAl2011}, dissipation-based analysis and design \cite{Ebenbauer2006}, controller synthesis \cite{jarvisEtAl2003}, reachability analysis \cite{Yin2021}, state observer synthesis~\cite{tan2006}, as well as synthesis of  barrier functions \cite{prajnaSafetyVerificationHybrid2004}, control Lyapunov functions \cite{tanSearchingControlLyapunov}, and control barrier functions \cite{amesEtAl2019}. 

Existing SOS toolboxes\deleted[id=TC]{, i.e., SOS parser,} such as \cite{Papachristodoulou2021SOSTOOLSMATLAB, Tobenkin2013SPOTless:Optimization, loefberg2009, Seiler2010a} are primarily designed to deal with convex problems. 
While convex SOS problems can be reduced to and solved as semidefinite programs (SDP) using state-of-the-art conic solvers \cite{Andersen2000TheAlgorithm, sedumi, ODonoghue2016},
\replaced[id=TC]{many applications}{most of the examples above} result in nonlinear (nonconvex) SOS problems. Solving nonconvex problems, however, using iterative \cite{chakrabortyEtAl2011} or sequential \cite{Cunis2023acc} methods requires the repeated evaluation of one or more parametrized convex subproblems. \added[id=TC]{When using the currently available toolboxes}, this necessitates the repeated re-initialization and parsing of the \added[id=JO2]{SOS problem, which can be time-consuming. Here, we refer \added[id=TC]{to} parsing as the transcription process (including polynomial operations) \added[id=TC]{of the sum-of-squares problem into an} SDP and \added[id=TC]{retrieval of the polynomial solution}. Recently, a new data structure was introduced \cite{Jagt2022EfficientSOSTOOLS} that significantly improves the computation time of \added[id=TC]{affine polynomial operations}. In case of nonconvex problems, \added[id=TC]{however,} the overall parsing time \added[id=TC]{is still high}.}

\begin{figure}[t]
    \center
    
    \includegraphics[width=\linewidth]{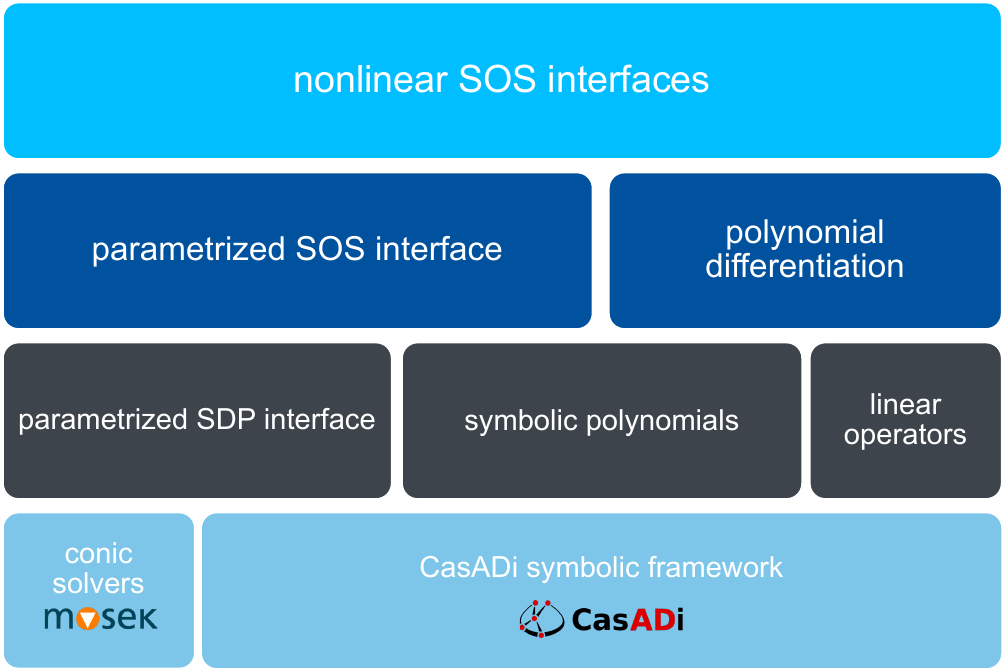}
    
    \caption{Overview of the Ca\textgreek{Σ}oS nonlinear sum-of-squares optimization suite and its main dependencies.}
    \label{fig:building-blocks}
\end{figure}

To improve the computation efficiency of nonlinear SOS \replaced[id=JO2]{\added[id=TC]{programs}}{solvers}, this paper introduces the nonlinear sum-of-squares optimization suite Ca\textgreek{Σ}oS. Ca\textgreek{Σ}oS provides a MATLAB-based framework for nonlinear symbolic polynomial expressions \deleted[id=TC]{which is} inspired by and built upon the symbolic framework of CasADi \cite{AnderssonCasADi-AControl} for automatic differentiation and nonlinear optimization.
Within this novel framework (Fig.~\ref{fig:building-blocks}), Ca\textgreek{Σ}oS offers interfaces for quasiconvex and nonlinear SOS problems, where the underlying parameterized (convex) SOS subproblem and its SDP representation are precomputed, thus avoiding the repeated re-initialization and parsing that slows down existing toolboxes. 
To that extent, Ca\textgreek{Σ}oS implements efficient polynomial data types, nonlinear polynomial operations, and linear operators\footnote{Linear operators realize, e.g., polynomial differentiation.} supported by a new concept of monomial sparsity. Additionally, Ca\textgreek{Σ}oS supports extensive conic vector, matrix, and polynomial constraints, including polynomial and exponential vector cones,\footnote{Only if supported by the conic solver.} diagonally dominant (DD) and scaled-diagonally dominant matrices (SDD), as well as DD-SOS and SDD-SOS relaxations \cite{ahmadiMajumdar2017}.


This paper is structured as follows: Section~\ref{sec: MathematicalBackground} provides the mathematical background on polynomial sum-of-squares optimization. Section~\ref{sec: Methodology} elaborates our methodology of a parametrized convex SOS problem, its SDP relaxation and illustrates the treatment of nonconvex SOS problems. Section~\ref{sec: ImplementationDetails} gives details about the implementation and provides a high-level comparison to existing SOS toolboxes. Section~\ref{sec: NumericalResults} provides numerical comparisons to other toolboxes.

\paragraph*{Notation}
$\mathbb S_n$ and $\mathbb S_n^+$ denote the sets of symmetric and positive semidefinite matrices, \added[id=JO]{respectively,} in $\mathbb R^{n \times n}$.

\section{Mathematical Background}
\label{sec: MathematicalBackground}
A polynomial up to degree $d \in \mathbb N_0$ is a linear combination
\begin{align}
    \label{eq:polynomial}
    \pi = \sum_{\| \alpha \|_1 \leq d} c_\alpha x^\alpha
\end{align}
where $x = (x_1, \ldots, x_\ell)$ is a tuple of $\ell \in \mathbb N$ indeterminate variables, $\{ c_\alpha \}_\alpha \subset \mathbb R$ are the coefficients of $\pi$, and $\alpha = (\alpha_1, \ldots, \alpha_\ell) \in \mathbb N_0^\ell$ are the degrees; we adopt the multi-index notation $x^\alpha := x_1^{\alpha_1} \cdots x_\ell^{\alpha_\ell}$. The set of polynomials in $x$ with real coefficients up to degree $d$ is denoted by $\mathbb R_d[x]$. The set of monomials $\{x^\alpha\}_{\| \alpha \|_1 \leq d}$ forms a basis for $\mathbb R_d[x]$. The polynomial $\pi$ can hence be associated with its {\em coordinate vector} $[c_\alpha]_{\| \alpha \|_1 \leq d} \in \mathbb R^{n_d}$. The set of polynomials of arbitrary degree, $\mathbb R[x] = \bigcup_{d \in \mathbb N_0} \mathbb R_d[x]$, is closed under addition, (scalar) multiplication, and integration and differentiation with respect to $x$. 

\begin{definition}
    A polynomial $\pi \in \mathbb R_{2d}[x]$ is a {\em sum-of-squares polynomial} ($\pi \in \Sigma_{2d}[x]$) if and only if there exist $m \in \mathbb N$ and $\pi_1, \ldots, \pi_m \in \mathbb R_d[x]$ such that $\pi = \sum_{i=1}^m (\pi_i)^2$.
\end{definition}

The seminal result of \cite[Theorem~3.3]{Parrilo2003} states that a polynomial $\pi \in \mathbb R_{2d}[x]$ is a sum-of-squares polynomial if and only if 
\begin{align}
    \label{eq:polynomial-gram}
    \pi = \sum_{{\|\alpha\|_1, \|\beta\|_1 \leq d}} Q_{\alpha,\beta} (x^\alpha x^\beta)
\end{align}
for some matrix $Q \in \mathbb S_{n_d}^+$, although the matrix $Q$ is usually non-unique. The set $\{ x^\alpha x^\beta \}_{\| \alpha \|_1, \| \beta \|_1 \leq d}$ is called the {\em Gram basis} of $\pi$. Both the vector basis and the Gram basis define linear mappings $z_d: \mathbb R^{n_d} \to \mathbb R_d[x]$ and $Z_d: \mathbb S_{n_d} \to \mathbb R_{2d}[x]$, respectively.
Note that $z_d$ is invertible whereas $Z_d$ is not.


The dual space to polynomials, 
the space $\mathbb R_d[x]^*$ of linear forms $\phi: \mathbb R_d[x] \to \mathbb R$,
is often associated with nonnegative (Borel) measures $\mu: \mathbb R^n \to \mathbb R_{\geq 0}$, viz. \cite{Lasserre2001}
\begin{align*}
    \langle \phi, \pi \rangle = \int \sum_\alpha c_\alpha x^\alpha \mathrm d \mu = \sum_\alpha c_\alpha y_\alpha
\end{align*}
where $y_\alpha = \int x^\alpha \mathrm d \mu$ are called {\em moments}; define $\mathbf y = \{y_\alpha\}_\alpha$. The dual cone $\Sigma_{2d}[x]^*$ is spanned by those Borel measures for which the moment matrix $M_d(\mathbf y) \in \mathbb S_{n_d}$ with $M_d(\mathbf y)_{\alpha,\beta} = y_{\alpha + \beta}$ is positive semidefinite.
{Again, we define $z_d^*: \mathbb R_d[x]^* \to \mathbb R_{n_d}$ and $Z_d^*: \mathbb R_{2d}[x]^* \to \mathbb S_{n_d}$.}

A (possibly nonlinear) sum-of-squares problem takes the form
\begin{align}
    \label{eq:sos-nonlinear}
    \min_{\xi \in \mathbb R_{2d}[x]^n} f(\xi) \quad \text{s.t. $\xi \in \Sigma_{2d}[x]^n$ and $g(\xi) \in \Sigma_{2d'}[x]^m$}
\end{align}
where $f: \mathbb R_{2d}[x]^n \to \mathbb R$ and $g: \mathbb R_{2d}[x]^n \to \mathbb R_{2d'}[x]^m$ are differentiable functions. Since the degrees $d, d' \in \mathbb N_0$ are fixed, $f$ and $g$ are mappings between finite vector spaces and the gradients $\nabla f(\cdot): \mathbb R_{2d}[x]^n \to \mathbb R$ and $\nabla g(\cdot): \mathbb R_{2d}[x]^n \to \mathbb R_{2d'}[x]^m$ are finite linear operators. Under a suitable constraint qualification, the optimal solution $\xi_0$ of \eqref{eq:sos-nonlinear} satisfies the Karush--Kuhn--Tucker (KKT) necessary conditions
\begin{align}
    \label{eq:kkt-nonlinear}
    \left\{
    \begin{alignedat}{2}
        &\nabla f(\xi_0) - \nabla g&(\xi_0)^* \lambda_g - \lambda_\xi = 0& \\
        &\langle \lambda_\xi, \xi_0 \rangle = 0, & \langle \lambda_g, g(\xi_0) \rangle = 0& \\
        &\xi_0 \in \Sigma_{2d}[x]^n, & g(\xi_0) \in \Sigma_{2d'}[x]^m&
    \end{alignedat}
    \right.
\end{align}
with dual variables $\lambda_\xi \in (\Sigma_{2d}[x]^n)^*$ and $\lambda_g \in (\Sigma_{2d'}[x]^m)^*$.


\section{Methodology}
\label{sec: Methodology}
To ease the solution of nonconvex optimization problems via iterative and sequential methods, we symbolically relax parametrized convex sum-of-squares problems to semidefinite programs. 

\subsection{Convex sum-of-squares problems}
To effectively solve \eqref{eq:sos-nonlinear}, we consider the parametrized convex sum-of-squares problem
\begin{multline}
    \label{eq:sos-convex}
    \min_{\xi \in \mathbb R_{2d}[x]^n} \tfrac{1}{2} \langle H(p) \xi, \xi \rangle + \langle \phi(p), \xi \rangle \\ \text{s.t. $A(p) \xi - \beta(p) \in \Sigma_{2d'}[x]^m$ and $\xi \in \Sigma_{2d}[x]^n$}
\end{multline}
where $H(\cdot): \mathbb R_{2d}[x]^n \to \mathbb (R_{2d}[x]^n)^*$, $\phi(\cdot): \mathbb R_{2d}[x]^n \to \mathbb R$ and $A(\cdot): \mathbb R_{2d}[x]^n \to \mathbb R_{2d'}[x]^m$ are linear operators with $H(\cdot)$ symmetric and positive definite\footnote{That is, $\langle H(\cdot) \xi, \zeta \rangle = \langle H(\cdot) \zeta, \xi \rangle$ for all $\xi, \zeta \in \mathbb R_d[x]^n$ and $\langle H(\cdot) \xi, \xi \rangle > 0$ for all $\xi \neq 0$.}, and $\beta(\cdot) \in \mathbb R_{2d'}[x]^m$ is a constant term, parametrized in $p$ (possibly a polynomial). As convex optimization, \eqref{eq:sos-convex} has the dual problem \cite{Dorn1960}
\begin{multline}
    \label{eq:sos-dual}
    \max_{\zeta \in \mathbb R_{2d}[x]^n, \lambda \in (\mathbb R_{2d'}[x]^m)^*} \langle \lambda, \beta(p) \rangle - \tfrac{1}{2} \langle H(p) \zeta, \zeta \rangle \\ 
    \text{s.t. $H(p) \zeta - A(p)^* \lambda + \phi(p) \in (\Sigma_{2d}[x]^n)^*$} \\ 
    \text{and $\lambda \in (\Sigma_{2d'}[x]^m)^*$} 
\end{multline}
where $A(\cdot)^*: (\mathbb R_{2d'}[x]^m)^* \to (\mathbb R_{2d}[x]^n)^*$ is the adjoint operator to $A(\cdot)$.
If \eqref{eq:sos-convex} and \eqref{eq:sos-dual} have a strongly dual solution $(\xi_0, \zeta_0, \lambda_0)$, then $\xi_0 = \zeta_0$ and $\lambda_\xi = H(p) \zeta_0 - A(p)^* \lambda_0 + \phi(p)$ and $\lambda_g = \lambda_0$ is a solution to the KKT system \eqref{eq:kkt-nonlinear} of \eqref{eq:sos-convex}.

\subsubsection{Relaxation to semidefinite programs}
Using the equivalence of sum-of-squares polynomials and positive semidefinite matrices, the problem in \eqref{eq:sos-convex} can be reduced to the parametrized semidefinite program
\begin{multline}
    \label{eq:sdp-primal}
    \min_{\substack{Q = (Q_1, \ldots, Q_n) \\ P = (P_1, \ldots, P_m)}} \tfrac{1}{2} \langle (Z_d^* \circ H(p) \circ Z_d) Q, Q \rangle + \langle (\phi(p) \circ Z_d), Q \rangle \\
    \text{s.t. $(z_{2d'}^{-1} \circ A(p) \circ Z_d) Q - (z_{2d'}^{-1} \circ Z_{d'}) P = z_{2d'}^{-1}(\beta(p))$} \\
    \text{and $Q_1, \ldots, Q_n \in \mathbb S_{n_d}^+$, $P_1, \ldots, P_m \in \mathbb S_{n_{d'}}^+$}
\end{multline}
with the affine constraint of \eqref{eq:sos-convex} rewritten into the so-called implicit or kernel form.\footnote{The explicit or image form yields a set of linear matrix inequalities (LMI) instead of PSD constraints; however, this form is more computationally demanding for larger problems and therefore not supported \cite{Parrilo2003}.}
Similarly, the dual problem \eqref{eq:sos-dual} becomes
\begin{multline}
    \label{eq:sdp-dual}
    \max_{\substack{\mathbf v = (v_1,\ldots,v_n) \subset \mathbb R^{n_{2d}} \\ \mathbf y = (y_1,\ldots,y_m) \subset \mathbb R^{n_{2d'}}}} \langle \mathbf y, z_{2d}^{-1}(\beta(p)) \rangle - \tfrac{1}{2} \langle (z_{2d}^* \circ H(p) \circ z_{2d}) \mathbf v, \mathbf v \rangle \\
    \text{s.t. $(Z_d^* \circ H(p) \circ z_{2d}) \mathbf v - (Z_d^* \circ A(p)^* \circ z_{2d'}^{-*}) \mathbf y$} \\
    \text{${} + Z_d^*(\phi(p)) \subset \mathbb S_{n_d}$ and $(Z_{d'}^* \circ z_{2d'}^{-*}) \mathbf y \subset \mathbb S_{n_{d'}}$}
\end{multline}
which can be seen as the dual to \eqref{eq:sdp-primal} if $\mathbf v$ and $z_{2d}$ are replaced by $R \in \mathbb S_{n_d}$ and $Z_d$, respectively. Note that the linear operator $Z_{d'}^* \circ z_{2d'}^{-*}$ corresponds to the moment matrix mapping $\mathbf y \mapsto M_{d'}(\mathbf y)$.

The problems in \eqref{eq:sdp-primal} and \eqref{eq:sdp-dual} are again convex and hence, if strong duality holds, 
we can retrieve the KKT solution of \eqref{eq:sos-convex} as
\begin{alignat*}{2}
    &\xi_0 = Z_d(Q_0), & \lambda_g = z_{2d'}^{-*}(\mathbf y_0),& \\
    &\lambda_\xi = H(p) \xi_0 - {} &A(p)^* \lambda_g + \phi(p)&
\end{alignat*}
using again the relationship between strongly dual and KKT solutions.

\subsubsection{Convex solver interfaces} 
Parametrized convex sum-of-squares problems such as \eqref{eq:sos-convex} are solved using the interface {\tt sossol} in Ca\textgreek{Σ}oS, where the problem is described symbolically, that is, without need to explicitly define $H(\cdot)$, $A(\cdot)$, $\phi(\cdot)$, or $\beta(\cdot)$. Conic problems with semidefinite matrix constraints and/or linear matrix inequalities are either solved using the symbolic interface {\tt sdpsol} or the lower-level interface {\tt conic}; these two can be understood as extensions of CasADi's {\tt qpsol} and {\tt conic} interfaces, respectively; and support the conic solvers {Mosek} \cite{Andersen2000TheAlgorithm}, {SeDuMi} \cite{sedumi}, and {SCS} \cite{ODonoghue2016}. All three solver interfaces support alternative cone constraints, including diagonally-dominant and scaled diagonally-dominant sum-of-squares polynomials ({\tt sossol}), diagonally-dominant and scaled diagonally-dominant matrices ({\tt sdpsol}) as well as, e.g., Lorentz (second-order or quadratic) and rotated Lorentz cone constraints ({\tt sdpsol} and {\tt conic}).

\subsection{Nonconvex sum-of-squares problems}
\label{subsec: NonConvSOSProb}
The parametrized problem \eqref{eq:sos-convex} builds the foundation for bisections as well as iterative and sequential methods for nonconvex sum-of-squares problems.

\subsubsection{Quasiconvex bisection}
A generalized sum-of-squares problem takes the form \cite{seilerBalas2010}
\begin{multline}
    \label{eq:sos-quasiconvex}
    \min_{t \in \mathbb R, \xi \in \mathbb R_{d}[x]^n} t \\
    \text{s.t. $t B(\xi) - A(\xi) \in \Sigma_{2d'}[x]^m$} \\
    \text{and $B(\xi) \in \Sigma_{2d'}[x]^m$, $C(\xi) \in \Sigma_{2d''}[x]^l$}
\end{multline}
with affine operators $A, B: \mathbb R_d[x]^n \to \mathbb R_{2d'}[x]^m$ and $C: \mathbb R_d[x]^n \to \mathbb R_{2d''}[x]^l$. It can be shown \cite[Theorem~2]{seilerBalas2010} that the problem in \eqref{eq:sos-quasiconvex} is equivalent to
\begin{align*}
    \min_{\xi \in \mathbb R_{d}[x]^n} f(\xi) \quad \text{s.t. $C(\xi) \in \Sigma_{d''}[x]^l$}
\end{align*}
where $f: \mathbb R_d[x]^n \to \mathbb R$ is quasiconvex, that is, $f$ has convex sublevel sets. It is common practice in quasiconvex optimization \cite{Agrawal2020} to solve \eqref{eq:sos-quasiconvex} as bisection over the feasibility problem
\begin{multline*}
    P(t): \exists \xi \in \mathbb R_d[x]^n \, \text{s.t. $t B(\xi) - A(\xi) \in \Sigma_{2d'}[x]^m$} \\
    \text{and $B(\xi) \in \Sigma_{2d'}[x]^m$, $C(\xi) \in \Sigma_{2d''}[x]^l$}
\end{multline*}
and we observe that $P(t)$ is a special case of \eqref{eq:sos-convex} parametrized in $p = t$. Quasiconvex sum-of-squares problems in Ca\textgreek{Σ}oS are defined symbolically using the {\tt qcsossol} interface.

\subsubsection{Coordinate descent (`V-s-iteration')}
If the sum-of-squares problem \eqref{eq:sos-nonlinear} is neither convex nor quasiconvex but bilinear, that is, any decision variable $\xi_1, \ldots, \xi_n$ respectively enters linearly, the most common approach to date is a coordinate descent. This approach has also been called `V-s-iteration' \cite{chakrabortyEtAl2011, chakrabortyEtAl2011ii}. Here, the bilinear problem is replaced by an alternating sequence of (quasi)-convex sum-of-squares problems
\begin{multline*}
    \xi_i^\dagger \in \argmin_{\xi_i \in \mathbb R_d[x]} f(\xi_1^\dagger, \ldots, \xi_{i-1}^\dagger, \xi_i, \xi_{i+1}^\dagger, \ldots, \xi_n^\dagger) \\
    \text{s.t. $g(\xi_1^\dagger, \ldots, \xi_{i-1}^\dagger, \xi_i, \xi_{i+1}^\dagger, \ldots, \xi_n^\dagger) \in \Sigma_{d'}[x]^m$} \\ 
    \text{and $\xi_i \in \Sigma_d[x]$}
\end{multline*}
for $i \in \{1, \ldots, n\}$,
where each (quasi)-convex problem is parametrized in the initial guesses or previous solutions $(\xi_1^\dagger, \ldots, \xi_{i-1}^\dagger, \xi_{i+1}^\dagger, \ldots, \xi_n^\dagger)$.
The benefit of this approach, apart from its apparent simplicity, is that the set of solutions is always feasible for the bilinear constraints. However, it also requires a feasible initial guess and convergence cannot be guaranteed \cite{chakrabortyEtAl2011}. Moreover, a suitable decomposition into (quasi)-convex sum-of-squares problems must be provided by the user. Therefore, Ca\textgreek{Σ}oS does not provide a designated interface for coordinate descent; but we note that the provided interfaces for parametrized (quasi)-convex optimization are very suitable for this type of strategy.

\subsubsection{Sequential convex programming}
\label{subsubsec: SeqSOS}
In \cite{Cunis2023acc}, we proposed to solve nonlinear sum-of-squares problems by a sequence of convex problems obtained from linearization of $f$ and $g$ at a previous solution. Generally speaking, such a sequential approach involves solving the parametrized convex sum-of-squares problem
\begin{multline*}
    \min_{\xi \in \mathbb R_{2d}[x]^n} \tfrac{1}{2} \langle H_k \xi - \xi_k, \xi - \xi_k \rangle + \langle \nabla f(\xi_k), \xi - \xi_k \rangle \\
    \text{s.t. $g(\xi_k) + \nabla g(\xi_k)(\xi - \xi_k) \in \Sigma_{2d'}[x]^m$} \\
    \text{and $\xi \in \Sigma_{2d}[x]^n$}
\end{multline*}
where $\xi_k \in \mathbb R_{2d}[x]^n$ is a previous solution (or initial guess) and $H_k: \mathbb R_{2d}[x]^n \to (\mathbb R_{2d}[x]^n)^*$ is an approximation for the second derivative of the Lagrangian in \eqref{eq:kkt-nonlinear}. This approach is also known as {\em quasi-Newton} method \cite{Izmailov2014} and often provides good convergence properties if the initial guess is carefully chosen, but typically requires some form of globalization strategy if not. While the full details of the sequential method are beyond the scope of this paper, Ca\textgreek{Σ}oS offers the interface {\tt nlsossol} for nonlinear sum-of-squares problems.

\section{Implementation Details}
\label{sec: ImplementationDetails}
We give details on the implementation of polynomial data types in Ca\textgreek{Σ}oS as well as its CasADi-like functions and parametrized solvers. At the end of this section, we compare core features of existing toolboxes to Ca\textgreek{Σ}oS.

\subsection{Overview}
Ca\textgreek{Σ}oS does not only offer interfaces for (convex, quasiconvex, and nonconvex) sum-of-squares optimization but, similar to most available toolboxes, provides its own implementation of polynomial data types. The core concept underpinning Ca\textgreek{Σ}oS' polynomial types is a strong separation between the indeterminate variables $x$ and the symbolic variables representing the polynomial coefficients $c_\alpha$ (the decision variables), while supporting nonlinear polynomial expressions.\footnote{In the sense of Ca\textgreek{Σ}oS, a nonlinear polynomial expression is an expression involving one or more symbolic polynomials (polynomials with a priori unknown coefficients) as well as constant polynomials (polynomials with given coefficients), of which the result is another polynomial, and which depends nonlinearly on the symbolic polynomials.}
This distinguishes the implementation of polynomial types in Ca\textgreek{Σ}oS from \added[id=TC]{that of} {\sc multipoly} \deleted[id=TC]{class} \added[id=TC]{used} by {\sc sosopt} and {\sc SOSTOOLS}, which treats all unknown coefficients as indeterminate variables; but it is also distinct from {\sc SOSTOOLS}' new {\tt dpvar} type \cite{Jagt2022EfficientSOSTOOLS} or the implementation in {\sc SPOTless}, which are limited to linear (affine) problem formulations.

The second concept distinguishing Ca\textgreek{Σ}oS from other sum-of-squares toolboxes is that of a multidimensional basis, which we call {\em monomial sparsity} in reference to matrix sparsity patterns in nonlinear optimization. The purpose of monomial sparsity patterns is twofold: First, it provides an indefinite basis for polynomial expressions as the user specifies objective and constraints of the sum-of-squares problems. Second, it serves as a definite (finite) basis for polynomial decision variables, and thus facilitates for the conversion between a polynomial and its coordinates. It is the second property of monomial sparsity patterns that allows us to define functions between polynomials including, as special case, parametrized sum-of-squares problems.

\subsection{Polynomial data types}
The monomial {\tt Sparsity} type in Ca\textgreek{Σ}oS maps the monomials of a vector or matrix-valued polynomial to its entries. As such, monomial sparsity can be understood as polynomials with binary coefficients: If the $(i,j)$-th entry of the matrix sparsity pattern $c_\alpha \in \mathbb B^{n \times m}$ belonging to the monomial $x^\alpha$ is `1', then the $(i,j)$-th entry of the matrix-valued polynomial $\pi \in \mathbb R_d[x]^{n \times m}$ has the corresponding monomial. Unlike in the previous section, where we made the blanket assumption that the vector-valued decision variable $\xi$ is of degree $d$ in all entries, the monomial sparsity type allows us to more accurately describe the basis of multidimensional polynomials, and hence more efficiently convert between a convex sum-of-squares optimization problem and its relaxation to a semidefinite program. Moreover, the monomial sparsity type forms the basis for linear operators (such as the gradients $\nabla f$ and $\nabla g$) and functions between polynomials. 

\begin{figure}
    \center
    
    \includegraphics[width=\linewidth]{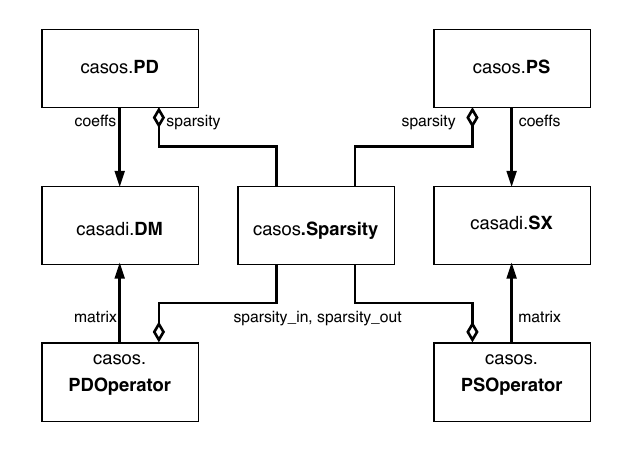}

    \caption{Relationship between monomial sparsity, polynomial data types, and linear operators in Ca\textgreek{Σ}oS as well as CasADi's matrix types.}
    \label{fig:polynomial-types}
\end{figure}

Fig.~\ref{fig:polynomial-types} illustrates the used data types and their relationships in Ca\textgreek{Σ}oS. With the multidimensional basis of a polynomial described by monomial sparsity, a polynomial can now be represented by the nonzero coordinate vector, that is, by the nonzero entries of the coefficients. Moreover, a (possibly nonlinear) polynomial expression can be reduced to the monomial sparsity pattern of its result and expressions for its nonzero coordinates. To that extent, Ca\textgreek{Σ}oS offers two polynomial data types: {\tt PD} for polynomials with constant (double) coefficients and {\tt PS} for polynomials of which the coefficients are symbolic variables and/or symbolic expressions. These two types extend and use CasADi's {\tt DM} and {\tt SX} matrix types (Fig.~\ref{fig:polynomial-types}), respectively.

At last, for a full polynomial algebra we need linear operators between polynomials. Examples for linear operators between polynomials include the gradients $\nabla f$ and $\nabla g$, the evaluation operator $E_a: \mathbb R_d[x] \to \mathbb R$ for any $a \in \mathbb R^\ell$, or the polynomial derivative (with respect to the indeterminate variables) $\partial_x: \mathbb R_d[x] \to \mathbb R_{d-1}[x]^\ell$. While linear operators between (finite) vectors can simply be represented by (finite) matrices, an example being the Jacobian matrix for a gradient, additional care is needed in the case of polynomials. In Ca\textgreek{Σ}oS, a linear operator $A: \mathbb R_d[x]^n \to \mathbb R_{d'}[x]^m$ is (internally) represented by the input and output bases (as monomial sparsity patterns) and a matrix $M \in \mathbb R^{n_{\mathrm o} \times n_{\mathrm i}}$ mapping the $n_{\mathrm i}$ nonzero input coordinates onto the $n_{\mathrm o}$ nonzero output coordinates. Two types of operators  are implemented: {\tt PDOperator} (linear operators between {\tt PD} polynomials with constant matrix $M$ of type {\tt DM}) and {\tt PSOperator} (linear operators between {\tt PS} polynomials with possibly symbolic matrix $M$ of type {\tt SX}). Unlike the functions described next, linear operators in Ca\textgreek{Σ}oS support a reduced, matrix-like interface including composition, addition, and scalar multiplication. 

Let us briefly explain how the linear operator for $\nabla g$ comes together: Suppose here that $g(\xi) \in \mathbb R_{d'}[x]^m$ is a nonlinear polynomial expression in the symbolic polynomial $\xi \in \mathbb R_d[x]^n$, where we require that all nonzero coordinates of $\xi$ are purely symbolic variables. Denote the nonzero coordinate vectors of $\xi$ and $g$ by $X \in \mathbb R^{n_X}$ and $G(X) \in \mathbb R^{n_G}$, respectively, where the latter is a symbolic expression in the purely symbolic vector $X$. The linear operator $\nabla g(\xi): \mathbb R_d[x]^n \to \mathbb R_{d'}[x]^m$ now is defined by the partial derivative $\nabla_X G(X) \in \mathbb R^{n_G \times n_X}$, the monomial sparsity of $X$ as input basis, and the monomial sparsity of $G(X)$ as output basis.

\subsection{Functions \& solvers}
Functions allow for the evaluation of polynomial expressions for given values (constant or symbolic) of the symbolic polynomials in the expression. This kind of evaluation ought not to be confused with the evaluation of a polynomial itself in its indeterminate variables or the substitution of indeterminate variables in a polynomial expression. With the monomial sparsity concept established in the previous subsection, functions between polynomials can be reduced to functions between their nonzero coordinates, and the Ca\textgreek{Σ}oS' {\tt Function} type hence makes use of the respective type in CasADi. 

\begin{remark}
    Although polynomials certainly can be interpreted as functions in the indeterminate variables and we could possibly evaluate a polynomial by substituting all its indeterminate variables by real numbers, the polynomial data types in Ca\textgreek{Σ}oS should not be treated as functions since they are not designed for evaluation. Instead, polynomial data types can be converted to CasADi {\tt Function} objects which do allow for rapid evaluation at multiple points at once.
\end{remark}

Similar to CasADi, the semidefinite and sum-of-squares solvers provided by Ca\textgreek{Σ}oS also offer the interface of a function, where values of the parameter, initial guesses, and further arguments are mapped to the optimal solution and its dual variables (optimality and feasibility provided). Here, unlike any other sum-of-squares toolbox, 
we explicitly separate between {\em initializing} a problem symbolically (e.g., establishing the Gram basis and precomputing operators or gradients) and {\em evaluating} its optimal solution. This strategy is particularly advantageous for nonconvex solvers, where the underlying convex subproblem is repeatedly solved. The result is a cascade of solvers, where the various parameters are propagated through the {\tt Function} interface.

\subsection{Comparison to other toolboxes}
Existing MATLAB toolboxes for automatic conversion of convex SOS problems to SDP include  {\sc SOSTOOLS} \cite{Papachristodoulou2021SOSTOOLSMATLAB}, {\sc sosopt} \cite{Seiler2010a}, both of which build upon {{\sc multipoly} \cite{seilerMultipolyToolboxMultivariable}}, a separate toolbox for polynomial operations and manipulations, {\sc SPOTless} \cite{Tobenkin2013SPOTless:Optimization}, and YALMIP with its SOS module \cite{loefberg2009}. Furthermore, there exist SOS toolboxes in programming languages other than MATLAB, e.g., {JuMP.jl \cite{Lubin2023} with its SOS extension. In the following, we limit the scope to MATLAB-based toolboxes to allow for a numerical comparison in the next section.

{\sc SOSTOOLS}, {\sc sosopt}, and YALMIP offer interfaces to various state-of-the-art conic solvers. All three provide methods to potentially improve the numerical conditioning as well as simplification techniques such as Newton polytopes to reduce the problem size. Notably, {\sc sosopt} also offers the dedicated bisection algorithm {\tt gsosopt} to solve quasiconvex sum-of-squares problems \cite{seilerBalas2010}. On the other hand, {\sc SPOTless} provides interface to only a few conic solvers but allows for DD-SOS and SDD-SOS constraints, which can be relaxed to linear and second-order cone constraints and scale better to larger problems at the cost of conservatism compared to the SDP relaxation \cite{ahmadiMajumdar2017}. Simplification techniques for problem size reduction are also avialable, but no dedicated methods for quasiconvex or nonconvex problems.
{The polynomial data type of {\sc multipoly} (used by {\sc SOSTOOLS} and {\sc sosopt}) offer a range of nonlinear polynomial operations at the cost of increased parsing times. The {\tt msspoly} class of {\sc SPOTless} and {\sc SOSTOOLS}'s new data structure {\tt dpvar} \cite{Jagt2022EfficientSOSTOOLS} are limited to linearity-preserving polynomial operations but their efficient implementations significantly improve the parsing process of convex problems.

 \comment[id=TC]{Maybe this would be good to discuss the numerical results?}
The toolboxes in this list share that none is \added[id=TC]{able to handle general} nonconvex SOS problems\deleted[id=TC]{ \added[id=TC]{with only} {\sc sosopt} \added[id=TC]{supporting at least} quasiconvex problems}. \deleted{Hence, the user always has to implement a coordinate descent algorithm manually.} \added[id=TC]{For iterative methods, hence}, \deleted[id=TC]{the \added[id=JO2]{SOS parsing process needs to be done} {convex SOS solvers have to be set up} in each iteration and the solution of one step has to be inserted into the other. In other words,} the conversion steps (\added[id=TC]{from SOS} to SDP and back to polynomial) must be called repeatedly \added[id=TC]{and a significant amount of} time is spent \added[id=TC]{parsing} the \added[id=TC]{sub}problem(s). 
This is a major difference to Ca\textgreek{Σ}oS where the parameterized conic solver(s) are \added[id=TC]{initialized} once \added[id=TC]{and repeated parsing is avoided}. \deleted[id=JO2]{Hence, during (e.g. coordinate descent) the main computational burden comes from the underlying SDP solver.} 
In addition, Ca\textgreek{Σ}oS offers a powerful nonlinear symbolic polynomial framework, parametrized simplification techniques, support of DD-SOS and SDD-SOS polynomial constraints,\footnote{Neither simplification techniques nor DD-SOS and SDD-SOS cones are discussed in this paper.} interfaces to the most common conic solvers, and methods for quasiconvex and nonconvex SOS problems.


    


\section{Numerical Results}
\label{sec: NumericalResults}
In this section we provide a numerical comparison of Ca\textgreek{Σ}oS to other MATLAB-based SOS toolboxes.
\added[id=JO2]{We consider three nonconvex benchmark problems: a region-of-attraction (ROA) problem \cite{chakrabortyEtAl2011}, a backwards-reachable set estimation problem \cite{Yin2021}, both for the longitudinal motion of an aircraft and a slightly modified ROA estimation problem for an N-link pendulum, to compare the scalability.} \added[id=JO2]{We compare the toolboxes regarding solution quality, iterations and/or computation time. Coordinate descent and bisections are used to solve the resulting, nonconvex SOS problems. Only Ca\textgreek{Σ}oS and {\sc sosopt} provide interfaces dedicated to quasiconvex problems, hence bisection algorithms are manually implemented for the other toolboxes.}

\begin{remark}
    Ca\textgreek{Σ}oS is the only toolbox offering sequential convex algorithms to deal with general nonconvex problems (see Section \ref{subsubsec: SeqSOS}). Hence, for comparison we implement a coordinate descent in combination with the quasiconvex module. See \cite{Cunis2023acc} for a comparison between coordinate descent/bisection and the sequential convex algorithms.
\end{remark}

\added[id=JO2]{We use the default options of each toolbox, except for {\sc SOSTOOLS} where we enabled the simplification based on Newton polytopes as they are activated per default also in other toolboxes.
\begin{remark}
    We have slightly adapted {\sc sosopt}, SOSTOOLS and SPOTless by measuring solver wall-time, disabling display messages, and augmented solver outputs. See the supplementary material as described below \cite{darus-4499} for more information.
\end{remark}
}

The results in this paper were computed on a personal computer \added[id=JO2]{with} Windows 10 and MATLAB 2023b running on an AMD Ryzen 9 5950X 16-Core Processor with 3.40 GHz and 128 GB RAM. We use Mosek \cite{Andersen2000TheAlgorithm} as the underlying SDP solver. We also provide supplementary material regarding the implementation details and the actual implementation of the benchmarks \cite{darus-4499}. 
\added[id=JO2]{Next, we shortly describe the procedure to measure the involved computation times. Afterwards, we provide the description and evaluation of the three distinct benchmarks.}

\subsection*{Time Measurement Procedure}
We distinguish between solve time -- the actual time spent in Mosek\footnote{We measure all solver instances that occur during the bisection. At the end we sum those up to obtain the total solve time.} -- and the \added[id=TC]{parse} time. The latter measures the time to setup the SOS problems (e.g., SOS constraints), the relaxation to SDP, and retrieving polynomial solutions. For all toolboxes except Ca\textgreek{Σ}oS we measure the total time spent in the framework including setup of the problem and subtract the solver time to obtain the \replaced[id=JO2]{parse}{build} time. For the benchmarks, we solve the problems five times with each toolbox and calculate the mean value of the computation times.
Since the SOS solver instances are only created once in Ca\textgreek{Σ}oS, we can directly measure its \replaced[id=JO2]{parse}{build} time. Additionally, we measure the time spent in a solver call and subtract the low-level solver time, which we refer to as \textit{call-time}. For simplicity, this call time is added to the \replaced[id=JO2]{parse}{build} time of Ca\textgreek{Σ}oS. 

\subsection{Benchmark I: Region-of-attraction estimation}
\added[id=JO2]{ The first benchmark deals with the computation of an inner-estimate of the ROA for the longitudinal motion of an aircraft with four states. In \cite{chakrabortyEtAl2011} both the nonlinear polynomial dynamics and the nonconvex SOS program are derived in terms of a three-step V-s-iteration algorithm (coordinate descent with bisection). We neither provide the derivation of the problem nor the actual SOS optimization problem due to limitation in space. For details, see \cite{chakrabortyEtAl2011} or the provided supplementary material \cite{darus-4499}.} \added[id=JO2]{The goal is to find an inner estimate of the region of attraction that reads $\Omega_\gamma := \{x \in \mathbb{R}^n \mid V(x) \leq \gamma \}$, where $V: \mathbb R^n \rightarrow \mathbb R$ is a Lyapunov function for the nonlinear dynamics and $\gamma > 0$ is a stable level set. The level set
$\Omega_\beta := \{x \in \mathbb{R}^n \mid p(x) \leq \beta \}$ for a given ellipsoidal $p(x)$ and $\beta >0$ is used to improve the solution quality via a set inclusion, i.e., $\Omega_\beta \subseteq \Omega_\gamma$. Within the three steps, the first two steps try to maximize the $\gamma$ and $\beta$ level sets (via bisection) for a given $V$. The third step seeks for a new Lyapunov candidate for the nonlinear system. These steps are repeatedly performed either for a maximum number of iterations $N_{\text{iter}} $ (here $N_{\text{iter}}= 100$) or if the $\beta$ level set seems not to further increase, i.e., $|\beta^{(k)}-\beta^{(k-1)}| \leq \epsilon$, where $\epsilon = 10^{-3}$ and $k \leq N_{\text{iter}}$ is the current iterate. }

\subsection*{Evaluation and discussion}
\added[id=JO2]{We compare computation times for each toolbox/parser in Fig.~\ref{fig:solveTimeComp_GTM_roa}. A full breakdown can be found in Table \ref{tab: compROA}. Ca\textgreek{Σ}oS (A) solves the problem in less than 70$\%$ of the time of the closest competitor ({\sc SOSTOOLS} (B) using {\tt dpvar}). Especially the parsing time of about \SI{16}{s} in Ca\textgreek{Σ}oS marks a significant reduction in comparison to other toolboxes. All other toolboxes have to re-initialize a new solver instance, that is, do the parsing process in every iteration. 
Here, the bisections ($\gamma$-step and $\beta$-step) are particularly demanding as they involve many evaluations of the convex subproblems.
SOSTOOLS (B), with its efficient {\tt dpvar} data structure, and SPOTless (E) handle repeated evaluations still rather efficiently (see also the comparison to SOSTOOLS' older data structure {\tt pvar} (C)), but they cannot compete with Ca\textgreek{Σ}oS (A) in terms of computation time. In this example, YALMIP's (F) parse time considerably exceeds the others. While the actual low-level solver times are similar across all toolboxes, deviations in the solve time can be explained by different transcription processes. The number of iterations are almost the same except for YALMIP (F) and SPOTless (E), which needed the fewest number of iterations. }

\added[id=JO2]{ In the third and fourth column in Table~\ref{tab: compROA}, we compare the toolboxes in terms of solution quality (sublevel set sizes). Here, SPOTless (C) computes the largest $\gamma$ level set, whereas {\sc SOSTOOLS} (B, C) and {\sc sosopt} (D) have the largest $\beta$ level set. Ca\textgreek{Σ}oS' $\gamma$ level set is the third largest and its $\beta$ level set is only slightly smaller compared to its closest competitor {\sc SOSTOOLS} (B) for which computation time is significantly higher. }

\begin{figure}[h!]
    \setlength{\figH}{4.5cm}
    \setlength{\figW}{6.5cm}
    \centering
%
%
\definecolor{mycolor1}{rgb}{0.00000,0.44700,0.74100}%
\definecolor{mycolor2}{rgb}{0.85000,0.32500,0.09800}%
\begin{tikzpicture}

\begin{groupplot}[%
group style={
    group size=1 by 2,
    vertical sep=1pt,
    x descriptions at=edge bottom,
},
width=\figW,
at={(0\figW,0\figH)},
scale only axis,
xmin=0.5,
xmax=6.5,
xtick={1,2,3,4,5,6},
xticklabels={{A},{B},{C},{D},{E},{F}},
xticklabel style={rotate=0},
axis background/.style={fill=white},
xmajorgrids,
ymajorgrids,
]

\nextgroupplot[
height=0.316\figH,
ymin=5600,
ymax=5850,
ytick={5600,5800},
bar width=0.4,
legend style={legend pos=north west, legend cell align=left, align=left, draw=white!15!black}
]
\addplot[ybar stacked, fill=mycolor1, draw=black, area legend] table[row sep=crcr] {%
6	5652.50273912\\
};
\addlegendentry{Parse time}

\addplot[ybar stacked, fill=mycolor2, draw=black, area legend] table[row sep=crcr] {%
6	96.3750544\\
};
\addlegendentry{Solve time}

\addplot[forget plot, color=white!15!black] table[row sep=crcr] {%
0.5	0\\
6.5	0\\
};

\nextgroupplot[
height=0.631\figH,
ymin=0,
ymax=500,
bar width=0.4,
ylabel style={font=\color{white!15!black}, align=right, text width=0.6\figH},
ylabel={Time (\si{\second})},
]
\addplot[ybar stacked, fill=mycolor1, draw=black, area legend] table[row sep=crcr] {%
1	15.89902714\\
2	62.31273914\\
3	171.31284308\\
4	308.58644658\\
5	82.8075289\\
6	5652.50273912\\
};
\addplot[forget plot, color=white!15!black] table[row sep=crcr] {%
0.5	0\\
6.5	0\\
};

\addplot[ybar stacked, fill=mycolor2, draw=black, area legend] table[row sep=crcr] {%
1	96.9693784\\
2	101.10981296\\
3	101.6777396\\
4	101.614405\\
5	99.2206\\
6	96.3750544\\
};
\addplot[forget plot, color=white!15!black] table[row sep=crcr] {%
0.5	0\\
6.5	0\\
};

\end{groupplot}
\end{tikzpicture}%
    \caption{Mean solve and parsing times for the ROA estimation problem using A:~Ca\textgreek{Σ}oS, B:~{\sc SOSTOOLS} (dpvar), C:~{\sc SOSTOOLS} (pvar), D:~{\sc sosopt}, E:~{\sc SPOTless}, F:~{\sc YALMIP}.}
    \label{fig:solveTimeComp_GTM_roa}
\end{figure}
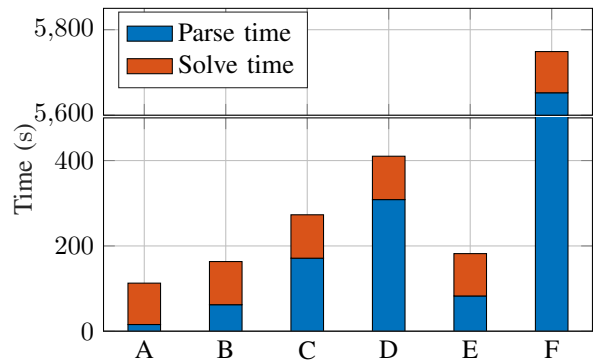

\begin{table}[htb!]
    \centering
    \caption{Comparison of the first benchmark regarding solution quality and mean computation times in seconds.}
    \begin{tabular}{lcccccc}
    \hline
    \hline
        Parser & Iterations & $\gamma^*$ & $\beta^*$ & Parsing  & Mosek   & Total   \\
    \hline             
        A & 55 & 169.9 & 274.9 & 15.9 & 96.97 & 112.87 \\
        B & 56 & 158.33 & 275.39 & 62.31 & 101.11 & 163.42 \\
        C & 56 & 158.33 & 275.39 & 171.31 & 101.68 & 272.99 \\
        D & 56 & 158.33 & 275.39 & 308.59 & 101.61 & 410.2 \\
        E & 47 & 214.48 & 274.66 & 82.81 & 99.22 & 182.3 \\
        F & 52 & 186.89 & 274.9 & 5652.5 & 96.38 & 5748.9 \\
    \hline
    \hline
    \end{tabular}
    \label{tab: compROA}
\end{table}

\subsection{Benchmark II: Backwards reachable set estimation}
\added[id=JO2]{ The second benchmark problem deals with computing the backward reachable set for the longitudinal aircraft flight dynamics as presented in \cite[Section V.B]{Yin2021} for a time horizon of $T = \SI{3}{s} $ without disturbances. \newpage The primary goal is to synthesize a set in the form $\Omega^V_{t,\gamma} = \{x \in \mathbb R^n \mid V(x,t) \leq \gamma \}$ where $V:\mathbb R^n \times \mathbb R \rightarrow \mathbb R$ is a so called \textit{reachability storage function} and $\gamma > 0$ bounds the set. Additionally, the goal is to synthesize a control law (with respect to control constraints), i.e.,  $u(t) =k(t,x)$ where $k: \mathbb R \times \mathbb R^n \rightarrow \mathbb R^m$, that ensures the system’s trajectories remain within a pre-defined target tube $\Omega_{t,0}^r := \{x \in \mathbb R^n \mid r(t,x) \leq 0\}$. 
The problem uses an adapted version of the polynomial aircraft model found in \cite{chakrabortyEtAl2011}. For details the reader is referred to \cite[Section V.B]{Yin2021} or the supplementary material \cite{darus-4499}.} 
The resulting bilinear problem is solved via a two-step coordinate-descent algorithm with bisections \linebreak{\cite[Algorithm 1]{Yin2021}}. 
\added[id=JO2]{ Similar to the first benchmark, we compute the example with different toolboxes and compare the results. Due to the already high computational effort of YALMIP in the ROA benchmark, we excluded it  from this example. To limit the computational effort we also limit the maximum number of iterations of the coordinate-descent to $N_{\text{iter}} = 10$.} We compare computation time and solution quality, i.e., volume of the sublevel set\footnote{The volume is approximated using Monte Carlo sampling.} in Fig.~\ref{fig:solveTimeComp_GTM_reach} and Table~\ref{tab: compTimesReac}.
 
\subsection*{Evaluation and discussion}
Comparing the parse times, Ca\textgreek{Σ}oS (A) needs significantly less time compared to the other toolboxes. SPOTless (C), the closest competitor in terms of parse times, takes approximately 20 times longer than Ca\textgreek{Σ}oS (A) to parse the convex subproblems. Especially {\sc sosopt}'s (D) parse times by far exceeds the others. The time spent in the low level solver differs between the four toolboxes, which can primarily be traced back to the different parsing processes and simplification techniques. 
In this specific case, SOSTOOLS (B) outperforms the competition in terms of solve times (Mosek) followed by Ca\textgreek{Σ}oS (A). However, due to its significantly lower parsing times, Ca\textgreek{Σ}oS (A) clearly outperforms each toolbox in terms of total computation times, almost halving the times of existing toolboxes. 

Additionally, the estimated sublevel set volume of Ca\textgreek{Σ}oS (A) is larger than SOSTOOLS (B) and  {\sc sosopt}'s (D). We were not able to compute a volume for SPOTless (C) although iterations returned feasible solutions.

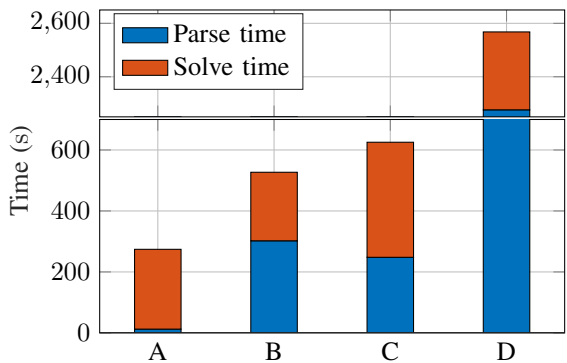
\begin{figure}[h!]
    \setlength{\figH}{4.5cm}
    \setlength{\figW}{6.5cm}
    \centering
%
%
\definecolor{mycolor1}{rgb}{0.00000,0.44700,0.74100}%
\definecolor{mycolor2}{rgb}{0.85000,0.32500,0.09800}%
\begin{tikzpicture}

\begin{groupplot}[%
group style={
    group size=1 by 2,
    vertical sep=1pt,
    x descriptions at=edge bottom,
},
width=0.951\figW,
at={(0\figW,0\figH)},
scale only axis,
xmin=0.5,
xmax=4.5,
xtick={1,2,3,4},
xticklabels={{A},{B},{C},{D}},
axis background/.style={fill=white},
xmajorgrids,
ymajorgrids,
]

\nextgroupplot[
height=0.316\figH,
ymin=2250,
ymax=2650,
ytick={2400,2600},
bar width=0.4,
legend style={legend pos=north west, legend cell align=left, align=left, draw=white!15!black}
]

\addplot[ybar stacked, fill=mycolor1, draw=black, area legend] table[row sep=crcr] {%
1	12.44589854\\
2	302.10359802\\
3	247.98530904\\
4	2276.26886188\\
};
\addplot[forget plot, color=white!15!black] table[row sep=crcr] {%
0.5	0\\
4.5	0\\
};
\addlegendentry{Parse time}

\addplot[ybar stacked, fill=mycolor2, draw=black, area legend] table[row sep=crcr] {%
1	261.75680654\\
2	224.83546736\\
3	377.4048\\
4	291.52853054\\
};
\addplot[forget plot, color=white!15!black] table[row sep=crcr] {%
0.5	0\\
4.5	0\\
};
\addlegendentry{Solve time}

\nextgroupplot[
height=0.631\figH,
ymin=0,
ymax=700,
bar width=0.4,
ylabel style={font=\color{white!15!black}, align=right, text width=0.6\figH},
ylabel={Time (\si{\second})},
]

\addplot[ybar stacked, fill=mycolor1, draw=black, area legend] table[row sep=crcr] {%
1	12.44589854\\
2	302.10359802\\
3	247.98530904\\
4	2276.26886188\\
};
\addplot[forget plot, color=white!15!black] table[row sep=crcr] {%
0.5	0\\
4.5	0\\
};

\addplot[ybar stacked, fill=mycolor2, draw=black, area legend] table[row sep=crcr] {%
1	261.75680654\\
2	224.83546736\\
3	377.4048\\
4	291.52853054\\
};
\addplot[forget plot, color=white!15!black] table[row sep=crcr] {%
0.5	0\\
4.5	0\\
};

\end{groupplot}
\end{tikzpicture}%
    \caption{Mean solve and parsing times for the reachable set estimation problem using A:~Ca\textgreek{Σ}oS, B:~{\sc SOSTOOLS} (dpvar), C:~{\sc SPOTless}, D:~{\sc sosopt}.}
    \label{fig:solveTimeComp_GTM_reach}
\end{figure}

\begin{table}[h!]
    \centering
    \caption{Comparison of the second benchmark regarding solution quality and mean computation times in seconds.}
    \begin{tabular}{lccccc}
    \hline
    \hline
        {Parser} &  {Volume}   &  {Parsing}  	& {Mosek} & {Total} \\
        \hline             
         A      &   $0.85(\pm 0.036)$   & 12.45 & 261.76   & 274.2 \\
         B 	   &   $0.79(\pm 0.035)$   & 302.1  & 224.84 & 526.94\\
         C      &   -                   & 247.99  & 377.4 & 625.39\\
         D 	&  $0.73(\pm 0.034)$       & 2276.27 & 291.53 & 2567.8\\
         \hline
         \hline
         \end{tabular}
    \label{tab: compTimesReac}
\end{table}

\subsection{Benchmark III: Region-of-attraction N-link pendulum}
The third and last benchmark compares the solve and parse times depending on the number of system states $n$. To this extent, we make use of an N-link pendulum (see e.g. \cite{MACHADO2016130} for the dynamics) which allows to easily scale the problem, i.e., one obtains $n = 2$N states.
 
Again, we consider the problem of computing the ROA with some small modification to the first benchmark. We fix $\gamma = 1$ and try to improve the size of the ROA by maximizing the size of the shape function. To limit the computation time, we set the maximum number of iterations for the coordinate descent to $N_{\text{iter}} = 20$. The parse and solver times depending on the system states are provided in Fig. \ref{fig:parseSolve_Nlink}.

 \added[id=JO2]{ While the solver times increase with similar magnitude for all toolboxes, the advantages of Ca\textgreek{Σ}oS (A) become more obvious as $n$ increases. Ca\textgreek{Σ}oS (A) performs best (overall solve time) due to its low parsing times, which remain consistently low across all problem sizes. Ca\textgreek{Σ}oS fundamental idea, of initializing (parsing) a problem symbolically once allows for efficient repeated numerical evaluations. This feature is clearly advantageous for solving nonconvex, potentially large-scale SOS problems. The major computational burden is now primarily shifted towards the solution method (e.g. bisection) and the actual low-level solver.}

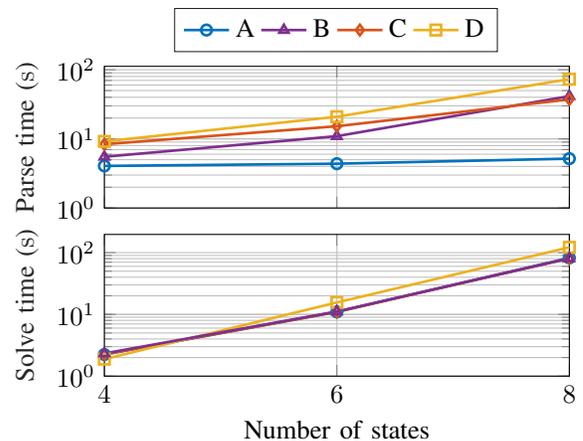
\begin{figure}[ht!]
    \setlength{\figH}{4.5cm}
    \setlength{\figW}{6.5cm}
    \centering
%
%
\definecolor{mycolor1}{rgb}{0.00000,0.44700,0.74100}%
\definecolor{mycolor2}{rgb}{0.85000,0.32500,0.09800}%
\definecolor{mycolor3}{rgb}{0.92900,0.69400,0.12500}%
\definecolor{mycolor4}{rgb}{0.49400,0.18400,0.55600}%
\begin{tikzpicture}

\begin{groupplot}[%
group style={
    group size=1 by 2,
    vertical sep=10pt,
    x descriptions at=edge bottom,
},
width=0.951\figW,
height=0.419\figH,
at={(0\figW,0.581\figH)},
scale only axis,
xmin=4,
xmax=8,
xtick={4,6,8},
ymode=log,
ymin=1,
yminorticks=true,
ylabel style={font=\color{white!15!black}},
axis background/.style={fill=white},
xmajorgrids,
ymajorgrids,
yminorgrids,
legend style={at={(0.15,1.1)}, anchor=south west, legend columns=-1, legend cell align=left, align=left, draw=white!15!black}
]

\nextgroupplot[ylabel={Parse time (\si{\second})},]

\addplot [color=mycolor1, line width=1.0pt, mark=o, mark options={solid, mycolor1}]
  table[row sep=crcr]{%
4	4.06219808\\
6	4.37455736\\
8	5.17930103999999\\
};
\addlegendentry{A}

\addplot [color=mycolor4, line width=1.0pt, mark=triangle, mark options={solid, mycolor4}]
  table[row sep=crcr]{%
4	5.51081052000001\\
6	10.94436936\\
8	41.4753463600001\\
};
\addlegendentry{B}

\addplot [color=mycolor2, line width=1pt, mark=diamond, mark options={solid, mycolor2}]
  table[row sep=crcr]{%
4	8.35603415999995\\
6	15.22871378\\
8	37.24157814\\
};
\addlegendentry{C}

\addplot [color=mycolor3, line width=1.0pt, mark=square, mark options={solid, mycolor3}]
  table[row sep=crcr]{%
4	9.19445733999998\\
6	20.91093816\\
8	73.0001474000001\\
};
\addlegendentry{D}

\nextgroupplot[ylabel={Solve time (\si{\second})}, xlabel={Number of states}]

\addplot [color=mycolor1, line width=1.0pt, mark=o, mark options={solid, mycolor1}]
  table[row sep=crcr]{%
4	2.23795406\\
6	10.98005066\\
8	82.1403566599999\\
};

\addplot [color=mycolor2, line width=1.0pt, mark=diamond, mark options={solid, mycolor2}]
  table[row sep=crcr]{%
4	2.25292266\\
6	10.96746088\\
8	80.9292600600001\\
};

\addplot [color=mycolor3, line width=1.0pt, mark=square, mark options={solid, mycolor3}]
  table[row sep=crcr]{%
4	1.88560000000003\\
6	15.604\\
8	121.8884\\
};

\addplot [color=mycolor4, line width=1.0pt, mark=triangle, mark options={solid, mycolor4}]
  table[row sep=crcr]{%
4	2.33070154\\
6	11.12738614\\
8	80.75376762\\
};

\end{groupplot}

\end{tikzpicture}%
    \caption{Mean solve and parsing times for the N-link pendulum ROA estimation problem using A:~Ca\textgreek{Σ}oS, B:~{\sc SOSTOOLS} (dpvar), C:~{\sc SPOTless}, D:~{\sc sosopt}.}
    \label{fig:parseSolve_Nlink}
\end{figure}

\addtolength{\textheight}{-4.6cm}   

\section{Conclusions}
{Ca\textgreek{Σ}oS is a MATLAB software suite specifically tailored to nonlinear sum-of-squares optimization. \added[id=TC]{By parsing repeatedly solved subproblems symbolically, combined with} efficient polynomial operations, Ca\textgreek{Σ}oS drastically reduces the overall computation time for \added[id=TC]{nonlinear SOS problems}. In benchmarks inspired from systems \added[id=TC]{and} control theory, \added[id=TC]{we} demonstrate that Ca\textgreek{Σ}oS is significantly faster \added[id=TC]{than} existing SOS toolboxes. \added[id=TC]{In particular,} the parsing time of Ca\textgreek{Σ}oS remains \added[id=TC]{virtually} constant even for a large number of states. \added[id=TC]{With Ca\textgreek{Σ}oS, we provide another step towards more tractable nonlinear sum-of-squares optimization for analysis and control of larger-scale systems.}}





\section*{Acknowledgment}
Fabian Geyer and  Renato Loureiro  contributed to the development of Ca\textgreek{Σ}oS.

\balance
\bibliographystyle{IEEEtran}
\bibliography{references,root}

\end{document}